# Smooth estimation of mean residual life under random censoring

Yogendra P. Chaubey[*,1] and Arusharka Sen[*,1]

*Concordia University*

**Abstract:** We propose here a smooth estimator of the mean residual life function based on randomly censored data. This is derived by smoothing the product-limit estimator using the Chaubey-Sen technique (Chaubey and Sen (1998)). The resulting estimator does not suffer from boundary bias as is the case with standard kernel smoothing. The asymptotic properties of the estimator are investigated. We establish strong uniform consistency and asymptotic normality. This complements the work of Chaubey and Sen (1999) which considered a similar estimation procedure in the case of complete data. It is seen that the properties are similar, though technically more difficult to prove, to those in the complete data case with appropriate modifications due to censoring.

## 1. Introduction

Let $T$ denote a non-negative random variable representing the lifetime of a subject or a component in the context of survival or reliability studies. Given that the component has survived up to time $t$, the expected remaining life is given by the so-called mean residual life (MRL) function

$$(1.1) \qquad m(t) = E(T - t | T > t) = \frac{\int_t^\infty S(u) du}{S(t)} I(S(t) > 0),$$

where $S(t)$ denotes the survival function (SF) of $T$.

Yang [22] was the first to propose estimating $m(t)$, by replacing the survival function in (1.1) with the empirical survival function (ESF). The reader may refer to Csörgo and Zitikis [6] and references therein for the vast literature on the properties of this estimator. There are several other articles that address the estimation of $m(t)$ from a random sample of censored lifetimes. For example, Yang [21], Hall and Wellner [12], Ghorai et al. [9] and Gill [11], among others, studied the properties of estimators of $m(t)$ based on randomly right-censored data, obtained generally by replacing $S(t)$ by its nonparametric estimator such as the Kaplan-Meier estimator or suitable modifications thereof. However, due to the inherent discontinuity of the MRL estimator described above, there is a natural interest in deriving smooth estimators of the MRL function when the random variable $T$ is assumed to admit a probability density function. In the complete-data case, Ruiz and Guillamòn [17] considered a kernel-smoothing-based estimator, and more recently, Abdous and Berred [1] have studied the properties of a local polynomial-based estimator of $m(t)$.

[*]Supported in part by a Discovery Grant from NSERC of Canada.
[1]Department of Mathematics and Statistics, Concordia University, Montreal, Quebec, Canada H3G 1M8, e-mail: chaubey@alcor.concordia.ca; asen@mathstat.concordia.ca

Chaubey and Sen [5] proposed an alternative MRL estimator in the complete-data case, based on a new smooth estimator of the survival function due to Chaubey and Sen [2]. Earlier, this smoothing technique was applied to hazard and cumulative hazard function estimation by Chaubey and Sen [3] in the complete-data case. Even though the method has been extended to density, survival, hazard and cumulative hazard function estimation for randomly censored data (Chaubey and Sen [4]), the MRL estimator has not been studied in the latter case.

Chaubey and Sen [5] realized that the weighting scheme proposed for smooth estimation of survival and density function was not appropriate for the smooth estimation of the MRL function and modified these weights. The purpose of the present paper is to adapt the above method for smooth estimation of mean residual life under the random censorship model. In this case we observe the censored lifetimes $Z_i = \min(T_i, C_i)$ and the indicator variable $\delta_i = I(T_i \leq C_i)$, where $C_1, C_2, \ldots$ is an independent sequence of non-negative independent random variables representing random censoring times. The Kaplan-Meier product-limit estimator (PLE) for the survival function is defined by

$$(1.2) \qquad \hat{S}_n(t) = \prod_{i=1}^{n} \left(1 - \frac{\delta_{[i:n]}}{n-i+1}\right)^{I_{[Z_{i:n} \leq t]}}$$

where $Z_{1:n} < \cdots < Z_{n:n}$ denote the ordered statistics corresponding to the observations $Z_i$, $i \leq n$ and $\delta_{[i:n]}$ is the value of $\delta$ corresponding to $Z_{i:n}$.

For the censored data, kernel-based smoothing has been adapted for density estimation as well as for hazard function estimation. For example, the kernel-based smooth estimation of the density function is given by (cf., for example, Mielniczuk [15])

$$(1.3) \qquad \hat{f}_n(x) = (\omega_n)^{-1} \int_0^\infty k\left(\frac{x-s}{\omega_n}\right) d\hat{F}_n(s),$$

where $\hat{F}_n(\cdot) = 1 - \hat{S}_n(\cdot)$, $\omega_n(>0)$, known as the *band-width*, is chosen so that $\omega_n \to 0$ but $n\omega_n \to \infty$, as $n \to \infty$; $k(.)$ is termed the *kernel function* and is typically assumed to be a symmetric pdf with zero mean and unit variance. On the other hand, hazard smoothing is based on the Nelson-Aalen estimator of the cumulative hazard function $H(x)$, namely,

$$H_n(x) = \sum_{i:Z_{i:n} \leq t} \frac{\delta_{(i:n)}}{n-i+1},$$

as given by

$$\begin{aligned} \hat{h}_n(t) &= (\omega_n)^{-1} \int_0^\infty k\left(\frac{t-s}{\omega_n}\right) d\hat{H}_n(s) \\ &= (\omega_n)^{-1} \sum_{i=1}^{n} \frac{\delta_{(i:n)}}{n-i+1} k\left(\frac{t-Z_{i:n}}{\omega_n}\right). \end{aligned}$$

A list of available results concerning the above smooth estimator and other variants is contained in the survey article by Padgett and McNichols [16].

Each of these methods can be used to derive a smooth estimator of mean residual life. The smoothing technique considered here differs from the above standard techniques in such a way that it does not suffer from problems faced by kernel smoothing for non-negative data, and derives the estimator of the MRL function as a plug-in estimator. The technique presented here could also be used in general to smooth the Nelson-Aalen estimator, but this is not considered here.

The form of the smooth estimator of the MRL function considered here is given in Section 2 and its asymptotic properties are investigated in Section 3. We also note that the results obtained here can be used to obtain the corresponding properties for the complete data case. In this process we have corrected the variance formula for the limiting asymptotic distribution of the resulting estimator.

## 2. A smooth estimator of the MRL function

As in Chaubey and Sen [4], the following representation of the product limit estimator may be used for computation. Let $M_n$ be the total number of failure points and let the ordered failure points be denoted by $Z^*_{i:n}$, $i = 1, 2, \ldots, M_n$, then we have

$$(2.1) \qquad \hat{S}_n(t) = \prod_{i \leq j} \left( \frac{n - k_i}{n - k_i + 1} \right), \ Z^*_{j:n} \leq t < Z^*_{j+1:n}, \ 0 \leq j \leq M_n$$

where $k_i$ is defined by $Z^*_{i:n} = Z_{k_i:n}$ ($Z^*_{0:n} = 0$, $Z^*_{M_n+1:n} = \infty$), i.e. $Z_{k_i:n}$ is the $k_i^{\text{th}}$ order statistic which corresponds to the $i^{\text{th}}$ failure.

For smooth estimators of the survival function and associated functionals in the random censorship case, Chaubey and Sen [4] consider the following modification of $S_n(t)$ as in Efron [7]:

$$(2.2) \qquad S_n(t) = \begin{cases} \hat{S}_n(t), & \text{for} \quad Z^*_{j:n} \leq t < Z^*_{j+1:n}, \ 0 \leq j \leq M_n \\ 0 & \text{for} \quad t > Z_{n:n} \end{cases}$$

where for $M_n < n$, we put $Z^*_{M_n+1:n} = Z_{n:n}$. With this modification the smooth estimator of the survival function proposed by Chaubey and Sen [4] is given by

$$(2.3) \qquad \tilde{S}_n(t) = \sum_{k=0}^{\infty} S_n(k/\lambda_n) p_k(\lambda_n x)$$

where

$$p_k(\mu) = e^{-\mu} \frac{\mu^k}{k!}, k = 0, 1, 2, \ldots,$$

$\lambda_n$ being a constant to be chosen suitably. The seemingly infinite series in (2.3) is in fact a finite sum with the last term being the $N^{\text{th}}$ term where $N = [\lambda_n Z_{n:n}]$, $[x]$ denoting the largest integer less than equal to $x$.

Chaubey and Sen [4] argue that the data-adaptive choice of $\lambda_n$ as given by $\lambda_n = n/Z_{n:n}$ may be appropriate for estimating the survival function. For studying the asymptotic properties of the resulting smooth estimator of the MRL function we will make a deterministic choice of $\lambda_n$, which may be obtained by cross-validation (see above).

Consistency of $\tilde{S}_n(t)$ can be proved using the same techniques as in Chaubey and Sen [4]. Towards this we extract the following theorem and corollary from Stute and Wang [20].

Let $F$ and $G$ denote the distribution functions of the survival time $T$ and censoring time $C$ and $H$ denote the common distribution of $Z_i, i = 1, 2, \ldots, n$, i.e. $1 - H = (1 - F)(1 - G)$. Set

$$\tau_H = \inf\{x : H(x) = 1\},$$

and write $F\{a\} = F(a) - F(a^-)$. Finally let $A$ denote the set of all atoms of $H$.

**Theorem 2.1** (Stute and Wang [20]). *Let $F_n = 1 - S_n$. Suppose that $F$ and $G$ have no jumps in common, and $\phi$ is $F$-integrable, then we have with probability one and in the mean,*

$$(2.4) \qquad \lim_{n \to \infty} \int \phi(x) F_n(dx) = \int_{x < \tau_H} \phi(x) F(dx) + 1_{(\tau_H \in A)} \phi(\tau_H) F\{\tau_H\}$$

*(where the second term should be zero if $\tau_H = \infty$).*

Two special cases of this result are of importance for further analysis and are given as corollaries below.

**Corollary 2.1.** *Suppose $F$ and $G$ have no jumps in common, then we have*

$$\sup_{t \leq \tau_H} |S_n(t) - S^*(t)| \to 0 \ \text{ with probability } 1,$$

*where*

$$S^*(t) = \begin{cases} S(t) & \text{if } t < \tau_H, \\ S(\tau_H^-) + 1_{\tau_H \in A} F\{\tau_H\}, & \text{if } t \geq \tau_H. \end{cases}$$

The uniform consistency of $S_n$ holds on $(0, \tau_H)$ if either $F\{\tau_H\} = 0$ or $F\{\tau_H\} > 0$ and $G\{\tau_H^-\} < 1$.

**Remark 2.1.** We will consider both $F$ and $G$ to be absolutely continuous hence, the uniform consistency for the product limit estimator will be assumed to hold on the support of $H$. In addition, if $\tau_H = \tau_F$ and $F$ is continuous (see Remark 3 of Stute and Wang [20]) then we have

$$(2.5) \qquad \lim_{n \to \infty} \int \phi(x) F_n(dx) = \int \phi(x) F(dx).$$

This implies the following corollary.

**Corollary 2.2.** *If $\tau_H = \tau_F$ and $F$ is continuous with $\int |x| F(dx) < \infty$, then with probability one we have*

$$(2.6) \qquad \lim_{n \to \infty} \int_{x > t} (x - t) F_n(dx) = \int_{x > t} (x - t) F(dx),$$

*where the above convergence is uniform on the interval $[0, \tau]$ for any $\tau < \tau_H$.*

**Remark 2.2.** If an addition $F$ is concentrated on non-negative values, the above implies that with probability one, $\lim_{n \to \infty} \int_{x > t} S_n(dx) = \int_{x > t} S(dx)$ uniformly on the interval $[0, \tau]$ for any $\tau < \tau_H$.

Now we are ready to establish the strong convergence of $\tilde{S}_n(t)$.

**Proposition 2.1.** *Let $H$ and $F$ have the same support and assume that $F$ is absolutely continuous. Consider a sequence of positive numbers $\lambda_n \to \infty$ as $n \to \infty$ then we have*

$$\sup_{t \leq \tau_H} |\tilde{S}_n(t) - S(t)| \to 0, \ \text{ as } n \to \infty.$$

*Proof.* The estimator considered here is based on the following theorem (see Lemma 1 in Feller [8], pp. 219):

Let $u(t)$ be a bounded and continuous function on $\mathbf{R}^+$. Then

$$(2.7) \qquad e^{-\lambda t} \sum_{k \geq 0} u(k/\lambda)(\lambda t)^k/k! \to u(t), \text{ as } \lambda \to \infty,$$

*uniformly in any finite interval* $\mathbf{J}$ *contained in* $\mathbf{R}^+$.

Using the above result, since $S(t)$ is bounded, continuous, nonnegative and non-increasing (on $\mathbf{R}^+$), we have

$$(2.8) \qquad S_n^*(t) = e^{-\lambda_n t} \sum_{k \geq 0} S(k/\lambda_n)(\lambda t)^k/k! \to S(t), \text{ as } \lambda_n \to \infty,$$

uniformly in any finite interval $\mathbf{J} \subseteq \mathbf{R}^+$.

Now we have

$$(2.9) \qquad \sup_{t \leq \tau_H} \{|\tilde{S}_n(t) - S_n^*(t)|\} \leq \max_{\{k:(k/\lambda_n) \leq \tau_H\}} |S_n(k/\lambda_n) - S(k/\lambda_n)| + S_n(\tau_H)\mathrm{Prob}\{N_n \geq [\tau_H \lambda_n]\},$$

where $N_n$ is a Poisson random variable with mean $t\lambda_n$. Both terms converge to zero with probability one, the first one due to Corollary 2.1 and the second one because $S_n(\tau_H)$ is bounded and $\mathrm{Prob}\{N_n \geq [\tau_H \lambda_n]\} \to 0$ as $\lambda_n \to \infty$. For analyzing the difference between $\tilde{S}_n(t)$ and $S(t)$, we see that

$$\sup_{t \leq \tau_H} |\tilde{S}_n(t) - S(t)| \leq \sup_{t \leq \tau_H} |\tilde{S}_n(t) - S_n^*(t)| + \sup_{t \leq \tau_H} |S_n^*(t) - S(t)|.$$

Using (2.8) and (2.9) we find that both terms on the right-hand side of the above inequality converge to zero and hence the theorem holds for $\tau_H < \infty$. For the case of $\tau_H = \infty$, the second term on the right-hand side of (2.9) is exactly zero and the theorem follows similarly. □

The estimator of the MRL function based on the Kaplan-Meier estimator as given by Yang [21] can be written as

$$m_n(t) = \frac{\int_t^\infty S_n(x)dx}{S_n(t)} I(t < Z_{M_n:n}).$$

We can similarly define the smooth estimator of the MRL function as

$$(2.10) \qquad \tilde{m}_n(t) = \frac{\int_{x>t}(x-t)\tilde{F}_n(dx)}{\tilde{S}_n(t)} = \frac{\int_t^\infty \tilde{S}_n(u)du}{\tilde{S}_n(t)},$$

where $\tilde{F}_n = 1 - \tilde{S}_n$.

The computational form of the above estimator is now provided in the following proposition as given in Chaubey and Sen [4].

**Proposition 2.2.** *The smooth estimator of the mean residual life as given by (2.10) can be represented as*

$$(2.11) \qquad \tilde{m}_n(t) = \frac{1}{\lambda_n} \frac{\sum_{k=0}^{N} \sum_{r=0}^{k} \frac{(t\lambda_n)^{k-r}}{(k-r)!} S_n(\frac{k}{\lambda_n})}{\sum_{k=0}^{N} \frac{(t\lambda_n)^k}{k!} S_n(\frac{k}{\lambda_n})},$$

*where* $N = [\lambda_n Z_{n:n}]$.

The representation given in (2.11) follows by direct integration and the fact that $S_n(x) = 0$ for $x > Z_{n:n}$.

## 3. Asymptotic properties of $\tilde{m}_n(\cdot)$

In this section we present the asymptotic properties of $\tilde{m}_n(\cdot)$. We establish its *strong uniform consistency* in Theorem 3.1 and asymptotic normality in Theorem 3.3.

**Theorem 3.1.** *Let $F$ and $G$ be absolutely continuous with common support and no common jump points such that $m(t) < \infty$ for all $t \in R^+$. Then for any compact interval $\mathcal{C} \subset \mathcal{R}^+$, we have for $\lambda_n \to \infty$ as $n \to \infty$*

$$\| \tilde{m}_n - m \|_{\mathcal{C}} = \sup_{t \in \mathcal{C}} |\tilde{m}_n(t) - m(t)| \to 0 \quad \text{a.s. as } n \to \infty. \tag{3.1}$$

*Proof.* If $F$ is compactly supported, i.e. $S(t) = 0$ for all $t > t_0 (< \infty)$, then $m(t)$ is not defined beyond $t_0$ and for all $t < t_0$, it is bounded. Thus we find that $m(t)$ is a Hadamard differentiable functional of $S(t)$. Hence, we can claim the same convergence property for $m_n(t)$, as that of $\tilde{S}_n(t)$ due to Proposition 2.1. Therefore, we confine ourselves to the case of infinite support, i.e. $\tau_H = \infty$.

First note that (2.11) can be written as

$$\tilde{m}_n(t) = \frac{(1/\lambda_n) \sum_{k=0}^{\infty} S_n(k/\lambda_n) P_k(t\lambda_n)}{\tilde{S}_n(t)}, \tag{3.2}$$

where $P_k(x) = \sum_{j \leq k} p_k(x)$. Writing the numerator in the above equation as $\sum_k \tilde{G}_n(k/\lambda_n) p_k(t\lambda_n)$, where $\tilde{G}_n(k/\lambda_n) = (1/\lambda_n) \sum_{j \leq k} S_n(j/\lambda_n)$, we note that

$$\tilde{G}_n\left(\frac{k}{\lambda_n}\right) - \frac{1}{\lambda_n} S_n\left(\frac{k}{\lambda_n}\right) \leq G_n^0\left(\frac{k}{\lambda_n}\right) \leq \tilde{G}_n\left(\frac{k}{\lambda_n}\right)$$

where

$$G_n^0(t) = \int_t^{\infty} S_n(x) dx. \tag{3.3}$$

Hence,

$$\sup_{t \in \mathcal{C}} |G_n^0\left(\frac{k}{\lambda_n}\right) - \tilde{G}_n\left(\frac{k}{\lambda_n}\right)| \to 0 \text{ a.s. as } n \to \infty, \tag{3.4}$$

where $\mathcal{C} = [0, \tau], \tau < \infty$. Further, since, $G_n^0(t)$ is non-decreasing and bounded *a.s.*, using Hille's theorem, we get

$$\sum_{k \geq 0} p_k(t\lambda_n) \tilde{G}_n(k/\lambda_n) \to G_n^0(t) \text{ a.s. as } n \to \infty \; \forall t \in \mathcal{C}. \tag{3.5}$$

Now using Corollary 2.2 (see (2.6)) we have

$$G_n^0(t) \to \int_t^{\infty} S(x) dx = m(t) \text{ a.s. as } n \to \infty. \tag{3.6}$$

Combining (3.5), (3.6) along with the fact that $\tilde{S}_n(t)^{-1} \to S(t)^{-1}$, uniformly in any compact interval we get the result in (3.1). $\square$

Before establishing the asymptotic normality of $\tilde{m}_n(t)$ we need the following notation. Let

$$H_0(z) = \mathrm{P}(Z \leq z, \delta = 0) = \int_0^z S(y) G(dy),$$

and
$$H_1(z) = \mathrm{P}(Z \leq z, \delta = 1) = \int_0^z (1 - G(y))F(dy).$$

Define
$$\gamma_0(x) = \exp\left\{\int_0^x \frac{H_0(dz)}{1 - H_1(z)}\right\},$$

$$\gamma_1(x) \equiv \gamma_\phi(x) = \frac{1}{1 - H(x)} \int_0^\infty I(x < w)\phi(w)\gamma_0(w)H_1(dw)$$

and
$$\gamma_2(x) \equiv \Gamma_\phi(x) = \int_0^\infty \int_0^\infty \frac{I(v < x, v < w)\phi(w)\gamma_0(w)}{(1 - H(v))^2} H_0(dv)H_1(dw).$$

Under the following assumptions,

(3.7)
$$\int \phi^2(x)\gamma_0^2(x)H_1(dx) < \infty$$

and

(3.8)
$$\int |\phi(x)|C^{1/2}(x)F^*(dx) < \infty,$$

where
$$C(x) = \int_0^x \frac{G(dy)}{(1 - G(y))(1 - H(y))},$$

Stute [19] proved the following CLT.

**Theorem 3.2.** *Under (3.7) and (3.8)*
$$\sqrt{n}\int \phi\, d(S_n - S^*) \longrightarrow N(0, \sigma_\phi^2),$$

*where*
$$\sigma_\phi^2 = \mathrm{Var}\{\phi(Z)\gamma_0(Z)\delta + \gamma_1(Z)(1 - \delta) - \gamma_2(Z)\}.$$

We show here that a similar central limit theorem holds for the smooth estimator $\tilde{m}_n(t)$. As mentioned in Remark 2.1, we suppose that $G$ is also continuous and assume the conditions (3.7) and (3.8) with $\phi(x) = x$, i.e.

(A.1)
$$\int_0^{\tau_H} \frac{x^2}{1 - G(x)} F(dx) < \infty,$$

and

(A.2)
$$\int_0^{\tau_H} x\sqrt{C(x)}F(dx) < \infty,$$

where $C(x)$ now becomes
$$C(x) = \int_0^x \frac{G(dy)}{(1 - H(y))(1 - G(y))}.$$

Note that the above two conditions imply

(B.1) A.1 holds, replacing $x^2$ in A.1 by $(\phi_t^A(x))^2$ or $(\phi_t^B(x))^2$ and A.2 holds replacing (the *first*) $x$ in A.2 by $\phi_t^A(x)$ or $\phi_t^B(x)$ for any $0 \leq t < \tau_H$, where $\phi_t^A(x) = (x-t)I(t > x)$, $\phi_t^B(x) = I(x > t)$.

(B.2) Given $\epsilon > 0$, there exists $M$, $t < M < \tau_H$, such that

$$\int_0^{\tau_H} \frac{x^2 I(x > M)}{1 - G(x)} F(dx) < \epsilon \text{ and } \int_0^{\tau_H} xI(x > M)\sqrt{C(x)}F(dx) < \epsilon.$$

The above equations are similar to (2.6) and (2.7) of Stute [19] (see p. 435). Further, note that

(B.3) for any $0 \leq a < \tau_H$, A.1 and A.2 hold for $\phi(x)I(x \leq a)$, i.e.

$$\int_0^{\tau_H} \frac{\phi^2(x)I(x \leq a)}{1 - G(x)} F(dx) < \infty,$$

and

$$\int_0^{\tau_H} |\phi(x)|I(x \leq a)\sqrt{C(x)}F(dx) < \infty,$$

whenever $\int \phi^2(x) F(dx) < \infty$. (see (2.3), p. 432 of Stute [19]).

**Theorem 3.3.** *Under assumptions A.1 and A.2 and if $\lambda_n \to \infty$, for any $t < \tau_H$*

$$\sqrt{n}(\tilde{m}_n(t) - m(t)) \longrightarrow N(0, \sigma^2(t)),$$

*where*

$$\sigma^2(t) = \text{ asymptotic variance of } \{\sqrt{n}(m_n(t) - m(t))\}.$$

(See Yang [21] or (3.11) below).

For further analysis we will use the following representation of the K-M estimator of the survival function as given in Stute [19]:

$$S_n(t) = \sum_{i=1}^n W_{in} I(Z_{i:n} > t),$$

where for $1 \leq i \leq n$,

$$W_{in} = \frac{\delta_{[i:n]}}{(n-i+1)} \prod_{j=1}^{i-1} \left[\frac{n-j}{n-j+1}\right]^{\delta_{[j:n]}}$$

is the mass attached to the $i$th order statistic $Z_{i:n}$ under $\hat{F}_n$. Therefore the smooth estimator of the survival function may be written as

$$\tilde{S}_n(t) = \sum_{i=1}^n W_{in} \phi_{nt}^B(Z_{i:n}),$$

where

$$\begin{aligned}\phi_{nt}^B(x) &= \sum_{k=0}^\infty p_k(t\lambda_n) I\left(x > \frac{k}{\lambda_n}\right) \\ &= \text{P}(\frac{N_n}{\lambda_n} < x),\end{aligned}$$

with $N_n \sim \text{Poisson}(t\lambda_n)$. Hence,

$$\begin{aligned}
\int_t^\infty \tilde{S}_n(y)dy &= \sum_{i=1}^n W_{in} \left[ \sum_{k=0}^\infty I(Z_{i:n} > \frac{k}{\lambda_n}) \int_t^\infty p_k(\lambda_n y)dy \right] \\
&= \sum_{i=1}^n W_{in} \left[ \sum_{k=0}^\infty I(Z_{i:n} > \frac{k}{\lambda_n}) \left\{ \frac{1}{\lambda_n} \sum_{r=0}^k p_r(t\lambda_n) \right\} \right] \\
&= \sum_{i=1}^n W_{in} \left[ \sum_{k=0}^\infty \left( \frac{[\lambda_n Z_{i:n}] - k + 1}{\lambda_n} \right) p_k(\lambda_n t) I(\lambda_n Z_{i:n} > k) \right] \\
&= \sum_{i=1}^n W_{in} \phi_{nt}^A(Z_{i:n}),
\end{aligned}$$

where

$$\begin{aligned}
\phi_{nt}^A(x) &= \sum_{k=0}^\infty \frac{[\lambda_n x] - k + 1}{\lambda_n} p_k(t\lambda_n) I\left(x > \frac{k}{\lambda_n}\right) \\
&= E_{N_n}\left[ \frac{[\lambda_n x] - N_n + 1}{\lambda_n} I(N_n < x\lambda_n) \right].
\end{aligned}$$

First note that, as $n \to \infty$,

$$\phi_{nt}^B(x) \to I(t < x) + \frac{1}{2} I(t = x) \text{ and } \phi_{nt}^A(x) \to (x - t)I(t < x).$$

Let us therefore consider,

$$(3.9) \quad S_n(t) = \sum_{i=1}^n W_{in} \phi_t^B(Z_{i:n})$$

$$(3.10) \quad \int_t^\infty S_n(y)dy = \sum_{i=1}^n W_{in} \phi_t^A(Z_{i:n}),$$

where

$$\phi_t^B(x) = I(x > t) \text{ and } \phi_t^A(x) = (x - t)I(x > t).$$

We have by Theorem 3.2

$$\sqrt{n}(m_n(t) - m(t))$$
$$= \sqrt{n}\left(\frac{\int_t^\infty S_n(y)dy}{S_n(t)} - m(t)\right)$$
$$= \sqrt{n}\frac{(\int_t^\infty S_n(y)dy - \int_t^\infty S(y)dy)}{S_n(t)} - \sqrt{n}m(t)\frac{(S_n(t) - S(t))}{S_n(t)}$$
$$= \sqrt{n}\frac{(\sum_{i=1}^n W_{in}\phi_t^A(Z_{i:n}) - E_F(\phi_t^A))}{S_n(t)}$$
$$\quad - \sqrt{n}m(t)\frac{(\sum_{i=1}^n W_{in}\phi_t^B(Z_{i:n}) - E_F(\phi_t^B))}{S_n(t)}$$
$$= \frac{1}{S(t)}\left[\frac{1}{\sqrt{n}}\sum_{i=1}^n\left(\frac{\delta_i\phi_t^A(Z_i)}{1-G(Z_i)} - E_F(\phi_t^A)\right)\right.$$
$$\quad + \frac{1}{\sqrt{n}}\sum_{i=1}^n\left((1-\delta_i)\gamma_{\phi_t^A}(Z_i) - \Gamma_{\phi_t^A}(Z_i)\right)$$
$$\quad - \frac{m(t)}{\sqrt{n}}\sum_{i=1}^n\left(\frac{\delta_i\phi_t^B(Z_i)}{1-G(Z_i)} - E_F(\phi_t^B)\right)$$
$$\left.\quad - \frac{m(t)}{\sqrt{n}}\sum_{i=1}^n\left((1-\delta_i)\gamma_{\phi_t^B}(Z_i) - \Gamma_{\phi_t^B}(Z_i)\right)\right]$$
$$(3.11) \qquad + o_P(1).$$

The claim in the above equation is justified by the assumptions in (3.7) and (3.8). In view of the above, to prove Theorem 3.2 it is enough to show that

$$(3.12) \qquad R_n(t) = \sqrt{n}\left[\frac{\int_t^\infty \tilde{S}_n(y)dy}{\tilde{S}_n(t)} - \frac{\int_t^\infty S_n(y)dy}{S_n(t)}\right] = o_P(1) \text{ as } n \to \infty.$$

To establish the above equation, we write

$$R_n(t) = \frac{\sqrt{n}[\int_t^\infty \tilde{S}_n(y)dy - \int_t^\infty S_n(y)dy]}{\tilde{S}_n(t)}$$
$$\quad - \sqrt{n}\frac{[\tilde{S}_n(t) - S_n(t)]}{\tilde{S}_n(t)}\left(\frac{\int_t^\infty S_n(y)dy}{S_n(t)}\right)$$
$$(3.13) \qquad = \frac{R_n^A(t)}{\tilde{S}_n(t)} - \left(\frac{\int_t^\infty S_n(y)dy}{S_n(t)}\right)\frac{R_n^B(t)}{\tilde{S}_n(t)}, \text{ say.}$$

To complete the proof of the assertion in (3.12), we establish the following lemmas.

**Lemma 3.1.** *Under A.1 and A.2,*

$$R_n^A(t) = o_P(1) \text{ as } n \to \infty.$$

**Lemma 3.2.** *Under A.1 and A.2,*

$$R_n^B(t) = o_P(1) \text{ as } n \to \infty.$$

*Proof of Lemma 3.1.* We have

$$
\begin{aligned}
R_n^A(t) &= \sqrt{n}\left[\int_t^\infty \tilde{S}_n(y)dy - \int_t^\infty S_n(y)dy\right] \\
&= \sum_{k=0}^\infty \sqrt{n}\int_t^\infty [S_n\left(\frac{k}{\lambda_n}\right) - S_n(y)]p_k(\lambda_n y)dy \\
&= \sum_{k=0}^\infty \int_t^\infty [\alpha_n\left(\frac{k}{\lambda_n}\right) - \alpha_n(y)]p_k(\lambda_n y)dy \\
&\quad + \sqrt{n}\sum_{k=0}^\infty \int_t^\infty [S\left(\frac{k}{\lambda_n}\right) - S(y)]p_k(\lambda_n y)dy,
\end{aligned}
$$

where

$$\alpha_n(y) = \sqrt{n}(S_n(y) - S(y)).$$

We may therefore write

$$
\begin{aligned}
R_n^A(t) &= \sum_{k=0}^\infty \int_t^M [\alpha_n\left(\frac{k}{\lambda_n}\right) - \alpha_n(y)]p_k(\lambda_n y)dy \\
&\quad + \sum_{k=0}^\infty \int_M^\infty \left[\alpha_n\left(\frac{k}{\lambda_n}\right) - \alpha_n(y)\right]p_k(\lambda_n y)dy \\
&\quad + \sqrt{n}\sum_{k=0}^\infty \int_t^\infty \left[S\left(\frac{k}{\lambda_n}\right) - S(y)\right]p_k(\lambda_n y)dy \\
&= r_{n,M}(t) + R_{n,M}(t) + A_n(t), \text{ say,}
\end{aligned}
\tag{3.14}
$$

where $M > M_1$ and $t < M_1 < \tau_H$ is chosen to satisfy B.2 for a given $\epsilon > 0$ and $A_n(t)$ is a deterministic function converging to zero, as shown in Chaubey and Sen [4]. Next, we have

$$
\begin{aligned}
r_{n,M}(t) &= \int_t^M \sum_{k:|\frac{k}{\lambda_n}-y|\leq h_n} \left[\alpha_n\left(\frac{k}{\lambda_n}\right) - \alpha_n(y)\right] p_k(\lambda_n y)dy \\
&\quad + \int_t^M \sum_{k:|\frac{k}{\lambda_n}-y|>h_n} \left[\alpha_n\left(\frac{k}{\lambda_n}\right) - \alpha_n(y)\right] p_k(\lambda_n y)dy \\
&= \int_t^M r_{n,M,1}^A(y)dy + \int_t^M r_{n,M,2}^A(y)dy, \text{ say.}
\end{aligned}
$$

By Corollary 3.1, p. 1316 of Deheuvels and Einmahl [10],

$$\sup_{t\leq y\leq M} |r_{n,M,1}^A(y)| = O(h_n[-\log h_n + \log\log n]) \tag{3.15}$$

provided $h_n \downarrow 0$, $\frac{nh_n}{\log n} \to \infty$ and $\frac{-\log h_n}{\log\log n} \to \infty$. Further, by (14), Inequality 1, p. 485 of Shorack and Wellner [18], we have

$$
\begin{aligned}
r_{n,M,2}^A(y) &\leq \sqrt{n}P[|\frac{N_n}{\lambda_n} - y| > h_n] \\
&\leq \sqrt{n}\exp\left(-\frac{\lambda_n h_n^2}{2y}\right).
\end{aligned}
\tag{3.16}
$$

It follows from (3.15) and (3.16) that

$$
\begin{aligned}
r_{n,M}^A(t) &\leq (M-t)O(h_n[-\log h_n + \log\log n]) \\
&\quad + (M-t)\sqrt{n}\exp\left(-\frac{\lambda_n h_n^2}{2M}\right) \text{ a.s.} \\
&\to 0 \text{ a.s. as } n \to \infty
\end{aligned}
$$

(3.17)

for any suitable choice of $h_n$.

Next consider $R_{n,M}^A(t)$. We have

$$
\begin{aligned}
R_{n,M}^A(t) &= \sum_{k=0}^{\infty}\int_M^{\infty}[\alpha_n\left(\frac{k}{\lambda_n}\right) - \alpha_n(y)]p_k(\lambda_n y)dy \\
&= \sqrt{n}\sum_{k=0}^{\infty}\int_M^{\infty}[S_n\left(\frac{k}{\lambda_n}\right) - S\left(\frac{k}{\lambda_n}\right)]p_k(\lambda_n y)dy \\
&\quad - \sqrt{n}\sum_{k=0}^{\infty}\int_M^{\infty}[S_n(y) - S(y)]p_k(\lambda_n y)dy.
\end{aligned}
$$

Now using the fact that

$$\int_M^{\infty} p_k(\lambda_n y)dy = \frac{1}{\lambda_n}\sum_{r=0}^{k}p_r(\lambda_n M)$$

in the first term on the RHS of the above equation, we have

$$
\begin{aligned}
R_{n,M}^A(t) &= \sqrt{n}\sum_{r=0}^{\infty}\left[\sum_{i=1}^{n}W_{in}\frac{([\lambda_n Z_{i:n}]+1-r)}{\lambda_n}I(Z_{i:n} > \frac{r}{\lambda_n})\right. \\
&\quad \left. - E_F\left(\frac{([\lambda_n X]+1-r)}{\lambda_n}I(X > \frac{r}{\lambda_n})\right)\right]p_r(\lambda_n M) \\
&\quad - \sqrt{n}\left[\sum_{i=1}^{n}W_{in}(Z_{i:n} - M)I(Z_{i:n} > M) - E_F((X-M)I(X > M))\right] \\
&= \sqrt{n}\sum_{r:(r/\lambda_n)\leq M_1}[\cdots] + \sqrt{n}\sum_{r:(r/\lambda_n)>M_1}[\cdots] \\
&\quad - \sqrt{n}\int \phi_M^A(x)(F_n(dx) - F(dx)) \\
&= R_{n,M,1}^A + R_{n,M,2}^A + R_{n,M,3}^A, \text{ say,}
\end{aligned}
$$

(3.18)

where $t < M_1 < M < \tau_H$ are as given after (3.14).

Note that the summands in $R_{n,M,1}^A$ and $R_{n,M,2}^A$ are of the form

$$\sqrt{n}\int \psi_{r,n}(dF_n - dF)p_r(\lambda_n M), \ r = 0, 1, 2, \ldots$$

where

$$
\begin{aligned}
\psi_{r,n}(x) &= \frac{([\lambda_n x]+1-r)}{\lambda_n}I(x > (r/\lambda_n)) \\
&= \left(x - \frac{r}{\lambda_N}\right)I\left(x > (r/\lambda_n)\right) + \frac{a_n(x)}{\lambda_n}I\left(x > (r/\lambda_n)\right),
\end{aligned}
$$

where $0 \leq a_n(x) < 1$. Furthermore, we have

$$\psi_{r,n}(x) = \phi_{n,(r/\lambda_n)}(x) + \frac{a_n(x)}{\lambda_n} I(x > (r/\lambda_n))$$

(3.19) $\quad \leq xI(x > M_1) + (1/\lambda_n)I(x > M_1) \text{ for } (r/\lambda_n) > M_1$

(3.20) $\quad \leq x + 1 \text{ for all } r = 0, 1, 2, \ldots$

Using the bounds in (3.19) and (3.20), and arguing as in the proof of Theorem 1.1, pp. 435–437 of Stute [19], we get by B.2 that

(3.21) $\quad \sqrt{n} \int \psi_{r,n}(dF_n - dF) = O_P(1) \text{ uniformly for } 0 \leq r \leq \lambda_n M_1$

(3.22) $\quad \sqrt{n} \int \psi_{r,n}(dF_n - dF) = O_P(\epsilon) \text{ uniformly for } r > \lambda_n M_1$

(3.23) $\quad \sqrt{n} \int \phi_M^A(dF_n - dF) = O_P(\epsilon)$

Thus

$$R_{n,M,1}^A = \sum_{r:(r/\lambda_n) \leq M_1} \sqrt{n}[\int \psi_{r,n}(dF_n - dF)] p_r(\lambda_n M)$$

$$= O_P(1) \text{Prob}\{\frac{N_{n,M}}{\lambda_n} \leq M_1\}, \text{ where } N_{n,M} \sim \text{Poisson}(\lambda_n M)$$

(3.24) $\quad = O_P(1)\text{Prob}\left\{\frac{N_{n,M} - \lambda_n M}{\sqrt{\lambda_n M}} \leq \sqrt{\lambda_n}\frac{M_1 - M}{M}\right\} \to 0 \text{ as } n \to \infty,$

since $M_1 - M < 0$ and $\lambda_n \to \infty$ as $n \to \infty$.

Also

$$R_{n,M,2}^A = \sum_{r:(r/\lambda_n) > M_1} \sqrt{n}[\int \psi_{r,n}(dF_n - dF)] p_r(\lambda_n M)$$

(3.25) $\quad = O_P(\epsilon),$

and $R_{n,M,3}^A = O_P(\epsilon)$ by (3.23).

Using (3.17), (3.18), (3.21)-(3.23), (3.24) and (3.25), Lemma 3.1 is established. □

*Proof of Lemma 3.2.*

$$R_n^B(t) = \sum_{k=0}^{\infty} \sqrt{n}\left[S_n(\frac{k}{\lambda_n}) - S_n(t)\right] p_k(\lambda_n t)$$

$$= \sum_{k=0}^{\infty} \left[\alpha_n(\frac{k}{\lambda_n}) - \alpha_n(t)\right] p_k(\lambda_n t)$$

$$+ \sum_{k=0}^{\infty} \sqrt{n}\left[S(\frac{k}{\lambda_n}) - S(t)\right] p_k(\lambda_n t),$$

so that negligibility of $R_n^B(t)$ follows by arguments similar to, but much easier than, those in the proof of Lemma 3.1. □

Now by (3.13), $R_n(t)$ is asymptotically negligible and hence the asymptotic normality of $\sqrt{n}[\tilde{m}_n(t) - m(t)]$ is now established and this also shows that the limiting distribution is the same as that of $\sqrt{n}[m_n(t) - m(t)]$.

**Remark 3.1.** In view of the above asymptotic result we note that a similar result holds in the complete data case, where $\delta \equiv 1$, $G(\cdot) \equiv 0$, $\Gamma_\phi(\cdot) \equiv 0$, so that by (3.11) the expression for the asymptotic variance simplifies to

$$
\begin{aligned}
&\sigma^2(t) \\
=\ & \text{var}\ ([(X-t)I(X>t) - m(t)I(X>t)]/S(t)) \\
=\ & (S^2(t))^{-1}[\text{var}((X-t)I(X>t)) + m^2(t)\text{var}(I(X>t)) \\
& -2m(t)\text{cov}((X-t)I(X>t),\ I(X>t))] \\
=\ & (S^2(t))^{-1}[E((X-t)^2 I(X>t)) - (\int_t^\infty S(y)dy)^2 + m^2(t)S(t)(1-S(t)) \\
& -2m(t)(1-S(t))\int_t^\infty S(y)dy] \\
=\ & (S^2(t))^{-1}[E((X-t)^2 I(X>t)) - m^2(t)S^2(t) + m^2(t)S(t)(1-S(t)) \\
& -2m^2(t)S(t)(1-S(t))],\ \text{since}\ \int_t^\infty S(y)dy = m(t)S(t), \\
=\ & (S^2(t))^{-1}[E((X-t)^2 I(X>t)) - m^2(t)S(t)] \\
=\ & \text{var}((X-t)|X>t)/S(t),
\end{aligned}
$$

which is the same as the one obtained by Yang [22]. Hence, for the complete data case we have

$$\sqrt{n}(\tilde{m}_n(t) - m(t)) \to N(0, \sigma^2(t)) \text{ as } n \to \infty.$$

This corrects the asymptotic distribution obtained in Chaubey and Sen [5].

## 4. Concluding remarks

We propose a smooth estimator of the mean residual life function under random censoring and investigate its asymptotic properties. To the best of our knowledge, no such estimator has been studied so far in the random-censoring case. We adapt the smoothing technique of Chaubey and Sen [2, 4], which is appropriate for non-negative data, especially in terms of avoiding boundary bias. We establish strong uniform consistency and an asymptotic weak representation of our estimator as an average of iid random variables. The latter, which automatically yields asymptotic normality, is established using the results of Stute [19] and Deheuvels and Einmahl [10]. The asymptotic variance of the estimator for uncensored data (Chaubey and Sen [5]) is derived as a special case.

**Acknowledgments.** The authors would like to thank a referee and Mervyn J. Silvapulle for their comments which have improved the readability of the original manuscript.